\documentclass{amsart}
\usepackage{graphicx}

\newtheorem{theorem}{Theorem}[section]

\newtheorem{lemma}[theorem]{Lemma}

\author{Sam Nelson}
\address{Department of Mathematics, Pomona College, 610 N. College
  Avenue, 
Claremont, CA 91711}
\email{knots@esotericka.org}
\author{Walter D.
  Neumann}
\address{Department of Mathematics, Barnard College,
  Columbia University, New York, NY 10027}
\email{neumann@math.columbia.edu}
\title{The 2-generalized knot group determines the knot} 
\dedicatory{To the memory of Xiao-Song Lin}
\begin{document}

\begin{abstract}Generalized knot groups $G_n(K)$ were introduced
  independently  by 
 Kelly (1991) and Wada (1992). We prove that $G_2(K)$ determines the
unoriented knot type and sketch a proof of the same for
$G_n(K)$ for $n>2$.
\end{abstract}
\maketitle
\section{The $2$--generalized knot group}

Generalized knot groups were introduced independently by Kelly
\cite{K} and Wada \cite{W}. Wada arrived at these group invariants of
knots by searching for homomorphisms of the braid group $B_n$ into
$\operatorname{Aut}(F_n)$, while Kelley's work was related to knot
racks or quandles \cite{FR,J} and Wirtinger-type presentations.

The Wirtinger presentation of a knot group expresses the group by
generators $x_1,\dots,x_k$ and relators $r_1,\dots,r_{k-1}$, in which
each $r_i$  has
the form $$x_j^{\pm1}x_ix_j^{\mp1}x_{i+1}^{-1}$$ for
some permutation $i\mapsto j$ of $\{1,\dots,k\}$ and map
$\{1,\dots,k\}\to\{\pm1\}$. The group $G_n(K)$ is defined by replacing
each $r_i$ by $$x_j^{\pm n}x_ix_j^{\mp  n}x_{i+1}^{-1}\,.$$ 
In particular, $G_1(K)$ is the usual knot group.

In \cite{T}, responding to a preprint of Xiao-Song Lin and the first
author \cite{LN}, Tuffley showed that $G_n(K)$ distinguishes the
square and granny knots.  $G_n(K)$ cannot distinguish a knot from its
mirror image. But $G_2(K)$ is, in fact, a complete unoriented knot
invariant.

\begin{theorem}
  The $2$--generalized knot group $G_2(K)$ determines the knot $K$ up
  to reflection.
\end{theorem}

We will assume $K$ is a non-trivial knot in the following proof, although it
is not essential. It is clear from the proof that the trivial knot is
the only knot with $G_2(K)=\mathbb Z$.

Wada described $G_n(K)$ as the fundamental group of the
space $M_n(K)$ obtained by gluing the boundary torus of the knot exterior to
another torus by the map $S^1\times S^1\to S^1\times S^1$ defined by
$f(z_1,z_2)=(z_1^n,z_2)$, where $z_1$ represents the meridian and
$z_2$ represents a longitude.  We will use this description. We will
call the glued-on torus the \emph{core torus}.

Note that $M_2(K)$ is a closed manifold: it can be described as the
result of gluing $Mb\times S^1$ into the knot exterior, where $Mb$
denotes the M\"obius band. It is clearly Haken, since its fundamental
group has a $\mathbb Z$ quotient, and it is irreducible and $\mathbb
P^2$--irreducible since its orientation cover is the double of the
knot exterior and hence irreducible. It therefore follows by Heil's
non-orientable extension  \cite{H} of Waldhausen's theorem
\cite{W} that $M_2(K)$ is determined by its fundamental group $G_2(K)$.

The core torus $T$ is the product $S^1\times S^1$ in $Mb\times
S^1\subset M_2(K)$,
where the first $S^1$ is the central circle of the M\"obius band. If
one cuts $M_2(K)$ along the core torus $T$ one recovers the knot
complement, from which the knot itself can be recovered by Gordon and
Luecke \cite{GL}. Thus the theorem follows from the following lemma.
\begin{lemma}
The core torus $T\subset M_2(K)$
is,  up to isotopy, the unique  one-sided torus in $M_2(K)$.
\end{lemma}
\begin{proof}
  We will use the ``geometric version'' of the JSJ decomposition. This
  is described, for instance, in \cite[Section 4]{NS}, but only for
  orientable manifolds, so we will discuss the non-orientable case
  briefly here. If the reader prefers to avoid JSJ for
  non-orientable manifolds (s)he can easily mirror our argument in the
  orientation cover of $M$.

  We restrict to the special case of an irreducible and $\mathbb
  P^2$--irreducible manifold $M$ whose boundary components are tori or
  Klein bottles.  The JSJ decomposition (geometric
  version\footnote{For manifolds covered by a torus bundle over the
    circle the JSJ decomposition described here is not
    necessarily the geometric one (the geometric JSJ decomposition is
    trivial in this case, overlooked in section 4 of
    \cite{NS}). It is an exercise to see that $M_2(K)$ is never of this form.})
  then decomposes $M$ along an embedded closed surface and is
  characterized by the first three of the following properties.
  \begin{enumerate}
  \item\label{1} The surface is a disjoint union of essential tori and Klein bottles,
  \item\label{2} $M$ is decomposed into simple (i.e., essential tori
    and annuli are boundary parallel) and Seifert fibered pieces, with
    no piece being an interval bundle over a torus or Klein bottle.
  \item The surface is minimal in the sense that it is, up to isotopy,
    a subset of any other surface with the above properties.
\item Any essential embedded torus or Klein bottle in $M$ can be
  isotoped to be a component of the JSJ surface, to lie in a neighborhood of
  one of its components, or to lie in a Seifert fibered
  piece of the decomposition.
\end{enumerate}
A short geometric proof of the existence and uniqueness of this
decomposition was given in \cite{NS} in the orientable case and can
easily be extended to the non-orientable case. Alternatively, one
decomposes the orientation cover of $M$ (or any other orientable
finite cover) and then uses the naturalness of the geometric JSJ
decomposition to descend to $M$. One of the features of the geometric
version of JSJ is that it lifts correctly in finite covers, and using
standard minimal surface technology one can isotop it to be preserved
by any finite group action.

Consider now the union of the JSJ surface $F$ for the knot exterior
and the core torus $T$. This is a surface that satisfies conditions
\eqref{1} and \eqref{2} so it contains the JSJ surface. It follows
easily that the JSJ surface is either $F$ or $F\cup T$. In the latter
case we know that any essential torus other than $T$ is isotopic into
the complement of $T$, hence embeds in $S^3$, and is thus two-sided,
so $T$ is the only one-sided torus. So assume the JSJ surface is
$F$. This only happens if the piece of the JSJ decomposition of the
knot exterior that contains the boundary is itself Seifert fibered, so $K$ is either
\begin{itemize}
\item  a $(p,q)$ torus knot for some $1<q<p$ (and $F$ is empty) or
\item a sum of $k>1$ prime knots.
\end{itemize}
In the first case $M$ is Seifert fibered over the orbifold $P(p,q)$
which is $\mathbb P^2$ with two orbifold points of degrees $p$ and $q$
respectively. In the second case the Seifert fibered component
containing $T$ fibers over an orbifold $Q(k)$ which is $k$-holed disk
with one boundary being a mirror boundary (the image of $T$).

Any essential surface in a Seifert fibered manifold is isotopic to a
vertical surface (union of fibers) or a transverse one (transverse to
all fibers). A transverse surface could only be closed in the first
case, but it is then hyperbolic since it covers the base orbifold
$P(p,q)$ which is hyperbolic. Thus in each case an essential embedded
torus must be vertical. In the first case, if it is one-sided it must
lie over an orientation reversing closed loop in $P(p,q)$, and there
is just one such loop up to isotopy (avoiding the orbifold points),
and it gives the torus $T$.  In the second case an essential torus
other than $T$ is the inverse image of a closed loop that does not
meet the mirror boundary or of a connected $1$-orbifold (i.e., an arc)
with both ends on the mirror boundary, and any such torus is
two-sided.
\end{proof}

\section{The $n$-generalized knot group}

The result holds also for $G_n(K)$ for $n>2$. Here is an outline
of the argument.

In this case $M_n(K)$ is not a manifold, so we cannot use
$3$--manifold JSJ. Instead we work directly with the group $G_n(K)$,
using the Scott-Swarup version of JSJ for groups \cite{SS}. For $n>3$
Scott-Swarup JSJ decomposes $M_n(K)$ as a graph of groups
corresponding to the JSJ decomposition of the knot exterior (in a
version close to classical JSJ rather than geometric JSJ), together
with an additional edge and vertex (of type ``$V_0$'' in the
terminology of SS-JSJ) as follows: the edge group and vertex group are
the peripheral subgroup of the knot exterior and $\pi_1(T)$. The edge
is characterized as the only edge of the graph of groups whose group
injects with finite index in a vertex group. The knot group can thus
be recovered as the fundamental group of the graph of groups which
results by removing this edge and its end vertex. Also, the peripheral
subgroup of the knot group is recovered as the edge group for this
edge. Finally, the knot is determined by knot group plus peripheral
subgroup by Gordon and Luecke \cite{GL}. For $n=3$ (and $2$) one can
use essentially the same argument, but there is an extra $V_0$--vertex
corresponding to the peripheral $\mathbb Z\times\mathbb Z$ of the knot
group, and the vertex for $\pi_1(T)$ is a $V_1$--vertex.

$G_n(K)$ is also defined for links, but is not a complete invariant.
Since it can be functorially derived from the rack (or quandle) of the
link, it cannot determine more than the rack determines (see
\cite{FR}). What $G_n(K)$ determines for a decomposable link is the
exteriors of the indecomposable sublinks, but since they are recovered
without knowledge of their orientation, one cannot reassemble the
whole link exterior. Moreover, since $G_n(K)$ (unlike the rack) does
not know the peripheral structure (i.e., the elements given by
meridians), it cannot always recover an indecomposable link,
since many links can share the same complement.

\subsection*{Acknowledgements} The second author acknowledges NSF
support under grant DMS-0456227. The authors thank Peter Scott, Colin
Rourke and Christopher Tuffley for useful correspondence. This paper
was born when the authors became acquainted on flight AA3394 (Baton
Rouge to Dallas 30/Mar/08), where they would not have met but for
Xiao-Song Lin, to whose memory they dedicate this paper.

\end{document}